\documentclass[11pt,draft]{article}
\usepackage{amsmath}
\usepackage{amssymb}
\usepackage{enumerate}
\usepackage{color}
\usepackage{xfrac}
\usepackage{bbm}
\usepackage[multiple]{footmisc}

\newtheorem{thm}{Theorem}[section]
\newtheorem{defn}[thm]{Definition}
\newtheorem{lemma}[thm]{Lemma}
\newtheorem{prop}[thm]{Proposition}
\newtheorem{rem}[thm]{Remark}
\newtheorem{cor}[thm]{Corollary}
\newtheorem{ex}[thm]{Example}

\newcommand{\bl}{\mathcal B_{\lambda}}
\newcommand{\tla}{\mathcal T_{\lambda,\alpha}}
\newcommand{\tFo}{{}_2F_1}
\newcommand{\eop}{\hfill$\square$}

\begin{document}
\title{Higher order Thorin-Bernstein Functions}

\author{Stamatis Koumandos and 
Henrik L. Pedersen\footnote{Research supported by grant DFF-1026-00267B from The Danish Council for Independent Research $|$ Natural Sciences} \footnote{Corresponding author}
}

\date{\today}
\maketitle

\begin{abstract}
We investigate subclasses of generalized Bernstein functions related to complete Bernstein and Thorin-Bernstein functions. Representations in terms of incomplete beta and gamma  as well as hypergeometric functions are presented. Several special cases and examples are discussed.
\end{abstract}
\noindent {\em \small 2020 Mathematics Subject Classification: Primary: 44A10. Secondary: 26A48, 33B15, 33C05.}

\noindent {\em \small Keywords: Laplace transform, completely monotonic function, Generalized Stieltjes function, Generalized Bernstein function, Thorin-Bernstein function} 


\newpage

\section{Introduction}

In this paper subclasses of the so-called generalized Bernstein functions of positive order are investigated. A non-negative function $f$ defined on $(0,\infty)$ having derivatives of all orders is a generalized Bernstein function of order $\lambda>0$ if $f'(x)x^{1-\lambda}$ is a completely monotonic function. This is equivalent to $f$ admitting an integral representation of the form 
\begin{equation}
f(x)=ax^{\lambda}+b+\int_0^{\infty}\gamma(\lambda,xt)\frac{d\mu(t)}{t^{\lambda}},
\label{eq:int-rep-bl}
\end{equation}
where $a$ and $b$ are non-negative numbers, and $\mu$ is a positive measure on $(0,\infty)$ making the integral converge for all $x>0$. 
Here, $\gamma$ denotes the incomplete gamma function 
$$
\gamma(\lambda,t)=\int_0^te^{-u}u^{\lambda -1}\, du.
$$
The convergence of the integral in \eqref{eq:int-rep-bl} is equivalent to  
$$
\int_0^{\infty}\frac{d\mu(t)}{(t+1)^{\lambda}}<\infty.
$$
The class of generalized Bernstein functions of order $\lambda$ is denoted by $\bl$ and was studied in \cite{KP1}. 

A function $g$ is called a generalized Stieltjes function of order $\lambda$ if there exist a positive measure $\nu$ on $[0,\infty)$ and a non-negative constant $c$ such that 
$$
g(x)=\int_0^{\infty}\frac{d\nu(t)}{(x+t)^{\lambda}}+c, 
$$
for $x>0$.  The class of these functions, denoted by $\mathcal S_{\lambda}$, can also be characterized in terms of the Laplace transform $\mathcal L$: $g\in \mathcal S_{\lambda}$ if and only if
$$
g(x)=\frac{1}{\Gamma(\lambda)}\mathcal L(t^{\lambda-1}\mathcal L(\nu)(t))(x)+c
=\frac{1}{\Gamma(\lambda)}\int_0^{\infty}e^{-xt}t^{\lambda-1}\mathcal L(\nu)(t)\, dt+c.
$$
This relation will be used throughout the paper. Several properties of generalized Stieltjes functions are given in \cite{KP2} and \cite{KP3}.
(For $\lambda=1$ the class is the class of Stieltjes functions, denoted by $\mathcal S$.)

The classes $\bl$ and $\mathcal S_{\lambda}$ are closely related: in fact, the map $\Phi:\bl\to \mathcal S_{\lambda}$ defined 
as
$$\Phi(f)(x)=x\mathcal L(f)(x)$$
is a bijection. See \cite[Theorem 3.1]{KP1}.

Let us also recall the definition of a completely monotonic function of order $\alpha$. A function $f:(0,\infty)\to [0,\infty)$ is completely monotonic of order $\alpha$ if $x^{\alpha}f(x)$ is completely monotonic. This class was introduced and characterized in \cite{KP0}. Let us remark, that $\alpha$ in this definition can be any real number. The class of completely monotonic functions of order $\alpha$ is denoted by $\mathcal C_{\alpha}$. 

The functions studied in this paper correspond to putting some extra conditions on the measure $\mu$ in the representation \eqref{eq:int-rep-bl}. 
\begin{defn}
\label{defn:fundamental}
 We say that a function $f:(0,\infty)\to \mathbb R$ is a higher order Thorin-Bernstein function, if there exist $\lambda>0$ and $\alpha<\lambda+1$ such that 
 \begin{equation}
 \label{eq:def-tb}
 f(x)=a x^{\lambda}+b+\int_0^{\infty}\gamma(\lambda,xt)\varphi(t)\, dt,
\end{equation}
where $a$ and $b$ are non-negative numbers, and $\varphi$ is a completely monotonic function of order $\alpha$.

For given numbers $\lambda$ and $\alpha$ the class of functions satisfying \eqref{eq:def-tb} is denoted by $\mathcal T_{\lambda,\alpha}$, and they are called $(\lambda,\alpha)$-Thorin-Bernstein functions.

\end{defn}
We notice that the integral in \eqref{eq:def-tb} converges exactly when $t^\lambda\varphi(t)/(t+1)^\lambda$ is integrable on $(0,\infty)$.

For given $\lambda>0$ it follows that $\mathcal T_{\lambda,\alpha_2}\subseteq \mathcal T_{\lambda,\alpha_1}$ when $\alpha_1\leq \alpha_2<\lambda +1$. When $\alpha=0$ the class $\mathcal T_{\lambda,\alpha}$ reduces to the class of generalized complete Bernstein functions of order $\lambda$, and when $\alpha=1$ to the class of generalized Thorin-Bernstein functions of order $\lambda$. These classes were studied in \cite{KP1}.

The higher order Thorin-Bernstein functions are closely related to the incomplete Beta function $B$, defined
 for $x\in [0,1)$, $a>0$, and $b\in \mathbb R$ as
$$
B(a,b;x)=\int_0^xt^{a-1}(1-t)^{b-1}\, dt.
$$
In fact, letting $\varphi(t)=t^{-\alpha}e^{-ct}=t^{-\alpha}\mathcal L(\epsilon_c)(t)$, for $c>0$, where $\epsilon_c$ denotes the point mass at $c$, it follows that
$$
\int_0^{\infty}\gamma(\lambda,xt)\varphi(t)\, dt=\Gamma(\lambda+1-\alpha)c^{\alpha-1}B\left(\lambda,1-\alpha;\frac{x}{x+c}\right),
$$
which is a special case of Proposition \ref{prop:incomplete-relation}. See also Example \ref{ex:fundamental}. Moreover, the incomplete Beta function is the main building block of the classes of higher order Thorin-Bernstein functions. See Corollary \ref{cor:A2}.

Bondesson introduced in \cite[p.~150]{B} the classes $\mathcal T_{1,\beta}$ and they also appear in \cite[Chapter 8]{S}. The higher order Thorin-Bernstein functions defined above play similar roles but in the context of generalized Bernstein functions of positive order $\lambda$. 

In this paper we shall relate these functions to, among other things, generalized Stieltjes functions and find representations in terms of hypergeometric functions.
Proposition \ref{prop:incomplete-relation} gives an equivalent integral representation of the functions $f\in \mathcal T_{\lambda,\alpha}$ in terms of the measure representing the completely monotonic function $t^{\alpha}\varphi(t)$ as a Laplace transform. This leads to a characterization in Theorem \ref{thm:A2} of the class $\mathcal T_{\lambda,\alpha}$ in terms of generalized Stieltjes functions of positive order: 
$f\in \tla \ \Leftrightarrow \ x^{1-\lambda}f'(x)\in \mathcal S_{\lambda+1-\alpha}$.

Based on the above result we prove that any generalized Bernstein function of order $\lambda>0$ is the pointwise limit of a sequence of functions from $\cup_{\alpha<\lambda +1}\mathcal T_{\lambda,\alpha}$. This result generalizes the corresponding result for Bernstein functions, obtained by Bondesson.

The motivation for studying these classes comes from specific examples including the incomplete Beta function. The incomplete Beta function plays, as the incomplete Gamma functions, an extensive role in probability and mathematical statistics. It also appears in e.g.\ Monte Carlo sampling in statistical mechanics (see \cite{kofte}).  For additional applications we refer to \cite{dlmf} and the references given therein.

In this paper all measures are supposed to be positive Radon measures.

\section{Fundamental results}

As mentioned above the classes $\mathcal T_{\lambda,\alpha}$ for fixed $\lambda>0$ are nested. Their intersection $\cap_{\alpha<\lambda+1}\mathcal T_{\lambda,\alpha}$ turns out to be the functions $\{ax^{\lambda}+b\, |\, a,b\geq 0\}$. This can be seen from Theorem \ref{thm:A2} as follows: if $f$ belongs to the intersection then $x^{1-\lambda}f'(x)$ belongs to $\cap_{\beta>0}\mathcal S_{\beta}$. This intersection, however, is equal to the constant functions (see \cite{Sokal}).

The class $\mathcal T_{\lambda,\alpha}$ is closed under pointwise limits:
\begin{prop}
\label{prop:easy-pointwise}
If $\{f_n\}$ is a sequence from $\tla$ that converges pointwise to $f:(0,\infty)\to \mathbb R$ then $f\in \tla$. 
\end{prop}
{\it Proof.}
Since $\tla$ is a subclass of $\bl$ and $\bl$ is closed under pointwise limits then $f\in \bl$ and furthermore, $f_n'\to f'$. See \cite[Proposition 2.4]{KP1}. This entails $x^{1-\lambda}f_n'(x)\to x^{1-\lambda}f'(x)$. From Theorem \ref{thm:A2} we also have $x^{1-\lambda}f_n'(x)\in \mathcal S_{\lambda+1-\alpha}$. Since the Stieltjes classes are also closed under pointwise limits (see \cite[Theorem 10]{KP}), $x^{1-\lambda}f'(x)\in \mathcal S_{\lambda+1-\alpha}$. Using again Theorem \ref{thm:A2} we see that $f\in \tla$.
\hfill $\square$

Since $\bl$ is closed under pointwise limits the limit of any pointwise convergent sequence of functions from  $\cup_{\alpha<\lambda +1}\mathcal T_{\lambda,\alpha}$ belongs to $\bl$. Conversely, we have the following proposition.
\begin{prop}
\label{prop:pointwise}
Let $f\in \bl$. Then there is a sequence from $\cup_{\alpha<\lambda +1}\mathcal T_{\lambda,\alpha}$ that converges pointwise to $f$.
\end{prop}

The proof of Proposition \ref{prop:pointwise} relies on the fact that any completely monotonic function is the pointwise limit of a sequence of generalized Stieltjes functions of positive order. See \cite{KP}. For the reader's convenience and in order to have a basis for proving Proposition \ref{prop:pointwise} we describe briefly a slightly different approach for proving this result from \cite{KP}. 

If $f$ is completely monotonic then by Bernstein's theorem, 
$$
f(x)=\int_0^{\infty}e^{-xt}\, d\sigma(t)
$$
for some positive measure $\sigma$ on $[0,\infty)$. Then take $\sigma_n=\sigma|_{[0,n]}$ and let 
$$
f_n(x)=\int_0^{\infty}e^{-xt}\, d\sigma_n(t)=
\int_0^{n}e^{-xt}\, d\sigma(t).
$$
These functions converge pointwise to $f$ (and even uniformly on $[r,\infty)$ for any $r>0$) and they are bounded by $\sigma([0,n])$.

Next define 
$$
h_k(s)=(1+s/k)^{-k}-e^{-s}, \ s\geq 0.
$$
It can be shown that $\{h_k\}$ is a non-negative and decreasing sequence of functions and furthermore that 
$$
\sup \{h_k(s)\, | \, s\geq 0\}\to 0
$$
as $k\to \infty$. (In fact, one may argue that $h_k$ has a unique maximum on $[0,\infty)$ and that the maximal value does not exceed $e/k$.)

Now, given a completely monotonic function $f$ we construct $\{f_n\}$ as above and then we choose, for any $n$ so large that $\sigma([0,n])>0$, $k_n$ such that  $$
\sup \{h_{k_n}(s)\, | \, s\geq 0\}\leq \frac{1}{n\sigma([0,n])}.
$$
Next we put $g_n(x)=\int_0^n(1+tx/k_n)^{-k_n}\, d\sigma(t)$ and notice that $g_n\in\mathcal S_{k_n}$ and also that $g_n-f_n$ converges uniformly to $0$, indeed:
$$
|f_n(x)-g_n(x)|\leq \int_0^nh_{k_n}(tx)\, d\sigma(t)\leq \frac{1}{n}.
$$
This gives
$$
g_n(x)=(g_n(x)-f_n(x))+f_n(x)\to f(x)
$$
as $n\to \infty$.

\noindent 
{\it Proof of Proposition \ref{prop:pointwise}.}
Let $f\in \bl$. Then $x^{1-\lambda}f'(x)$ is completely monotonic, and thus of the form 
$$
x^{1-\lambda}f'(x)=\int_0^{\infty}e^{-xt}\, d \sigma(t).
$$
Next let $f_n(x)=\int_0^{n}e^{-xt}\, d\sigma(t)$ and construct $g_n\in \mathcal S_{k_n}$ as above such that 
$$
\sup \{[f_n(x)-g_n(x)|\,|\, x\geq 0\}\to 0,\ \text{as}\ n\to \infty.
$$
Now let $x>0$ be given. Then $\int_0^xt^{\lambda-1}(f_n(t)-g_n(t))\, dt\to 0$, as $n\to \infty$. Indeed, the integrand tends pointwise to $0$ and 
$$
t^{\lambda-1}\sup_n\sup_{s\geq 0}|f_n(s)-g_n(s)|
$$
is an integrable majorant so that Lebesgue's theorem on dominated convergence can be applied.

Furthermore, 
$$
\int_0^xt^{\lambda-1}g_n(t)\, dt=\int_0^xt^{\lambda-1}(g_n(t)-f_n(t))\, dt+\int_0^xt^{\lambda-1}f_n(t)\, dt.
$$
The first term tends to $0$, and the second term tends, by monotone convergence, to 
$$
\int_0^xt^{\lambda-1}\int_0^{\infty}e^{-ts}\, d\sigma(s)\, dt=\int_0^xf'(t)\, dt=f(x)-f(0^+).
$$
Thus, defining
$$
F_n(x)=\int_0^xt^{\lambda-1}g_n(t)\, dt,
$$
then $F_n(x)\to f(x)-f(0^+)$ as $n\to \infty$ and $x^{1-\lambda}F_n'(x)=g_n(x) \in \mathcal S_{k_n}$. From Theorem \ref{thm:A2} it now follows that $F_n\in \mathcal T_{\lambda,\lambda+1-k_n}$. This completes the proof. \eop


By inspection of the construction preceeding the proof of Proposition \ref{prop:pointwise} one sees that any bounded completely monotonic function is the uniform limit of a sequence of bounded generalized Stieltjes functions. This observation can be used to obtain the next corollary, which refines \cite[Proposition 3.5]{KP1}).

\begin{cor}
\label{cor:finite-borel-measure}
     For any finite Borel measure $\mu$ on $[0,\infty)$ there exists a sequence of bounded generalized Bernstein functions $b_n$ of positive order and a sequence of non-negative numbers $c_n$  such that 
 $$
 c_n\epsilon_0+b_n'(x)\, dx\to \mu \ \text{weakly on}\ [0,\infty).
 $$
\end{cor}
{\it Proof.} Given $\mu$ we notice that $f=\mathcal L(\mu)$ is a bounded completely monotonic function. As remarked above there is a sequence $\{f_n\}$ where $f_n\in\mathcal S_{\lambda_n}$ is bounded such that $f_n\to f$ uniformly on $[0,\infty)$ as $n\to \infty$. The function $f_n$  can be represented as
$$
f_n(x)=c_n+\int_0^{\infty}\frac{d\mu_n(t)}{(x+t)^{\lambda_n}}=\mathcal L\left( \sigma_n\right)(x),
$$
where
$$
\sigma_n=c_n\epsilon_0+\frac{t^{\lambda_n-1}}{\Gamma(\lambda_n)}\mathcal L(\mu_n)(t)dt.
$$
Boundedness of $f_n$ means that $\sigma_n$ is a finite measure. Since $\mathcal L(\sigma_n)\to \mathcal L(\mu)$ pointwise on $[0,\infty)$ it follows that $\sigma_n\to \mu$ vaguely and $\sigma_n([0,\infty))\to \mu([0,\infty))$, and therefore $\sigma_n\to \mu$ weakly. 
Defining
$$
b_n(x)=\int_0^{x}\frac{t^{\lambda_n-1}}{\Gamma(\lambda_n)}\mathcal L(\mu_n)(t)dt,
$$
we see that $b_n\in \mathcal B_{\lambda_n}$, that $b_n$ bounded (indeed $b_n(x)\leq \sigma_n([0,\infty))$) and that $c_n\epsilon_0+b_n'(x)\, dx=\sigma_n$. This proves the assertion.
\eop

It should be noted that bounded generalized Bernstein functions have been characterized in  \cite[Proposition 3.9]{KP1}.
\begin{rem}
    If $f\in \bl$ is represented by the triple $(a,b,\mu)$ and if $f_n\in \mathcal T_{\lambda,\alpha_n}$, represented by the triple $(a_n,b_n,\mu_n)$, converges pointwise to $f$ then 
    $$
    x^{\lambda-\alpha_n}\mathcal L(\mu_n)(x)\, dx\rightarrow \mu \ \text{vaguely on}\ x>0.
    $$
    Indeed, we have according to \cite[Proposition 2.4]{KP1} that
    $$
    x^{1-\lambda}f_n'(x)\to x^{1-\lambda}f'(x)=\mathcal L(\lambda a\epsilon_0+\mu)(x),\ \text{as}\ n\to \infty.
    $$
    Also, using Theorem \ref{thm:A2}
    \begin{align*}
    x^{1-\lambda}f_n'(x)&=\lambda a_n+\Gamma(\lambda+1-\alpha_n)\int_0^{\infty}\frac{d\mu_n(s)}{(s+x)^{\lambda+1-\alpha_n}}\\
    &=\mathcal L(\lambda a_n \epsilon_0+s^{\lambda-\alpha_n}\mathcal L(\mu_n)(s))(x)
    \end{align*}
    and this gives 
    $$
    \lambda a_n \epsilon_0+s^{\lambda-\alpha_n}\mathcal L(\mu_n)(s)\,ds\rightarrow \lambda a\epsilon_0+\mu
    $$
    vaguely on $[0,\infty)$. See \cite[Proposition 9.5]{bf}. Restriction to the open half line establishes the result.
    \end{rem}

Next we develop the theory aiming, among other things, at the proof of Theorem \ref{thm:A2}.
The first lemma is a simple consequence of Fubini's theorem.
\begin{lemma}
 \label{lemma:F2}For any non-negative Borel measurable function $g$ and any positive measure $\mu$ we have
 $$
 \int_0^{\infty}g(t)\mathcal L(\mu)(t)\, dt=\int_0^{\infty}\mathcal L(g)(s)\, d\mu (s).
 $$
\end{lemma}

The next result is about Laplace transforms and convolution measures. For the reader's convenience let us mention that the convolution measure $\mu\ast\nu$ of two  measures $\mu$ and $\nu$ is given as the image measure of the product measure under the map $\tau(x,t)=x+t$.

\begin{lemma}\label{lemma:X} Let $\mu$ and $\nu$ be two measures on $[0,\infty)$ and let $\beta>0$. Then
$$
\int_0^{\infty}t^{\beta-1}\mathcal L(\mu)(t)\mathcal L(\nu)(t)\, dt=\Gamma(\beta)\int_0^{\infty}\frac{d(\mu\ast \nu)(s)}{s^{\beta}}.
$$
\end{lemma}
{\it Proof.} This follows from Lemma \ref{lemma:F2} and the fact that $\mathcal L(t^{\beta-1})(s)=\Gamma(\beta)/s^{\beta}$. 
\eop
\begin{rem}
 Because of positivity of the function and measure in Lemma \ref{lemma:X} interchanging the order of integration is permitted. It may of course happen that the integrals equal $\infty$. For values of $t$ near $0$, $\mathcal L(\mu)(t)\mathcal L(\nu)(t)$ is in general bounded from below by some positive constant. Hence, in order for the integral on the left hand side to be convergent, $\beta$ must be positive.
\end{rem}
\begin{cor}
\label{lemma:F3}
 Let $\lambda >0$ and $\alpha<\lambda +1$. For $f\in \bl$ and any measure $\mu$ on $[0,\infty)$ the following relation holds:
 \begin{align*}
 \int_0^{\infty}f(t)t^{-\alpha}\mathcal L (\mu)(t)\, dt&=\Gamma(\lambda+1-\alpha)\int_0^{\infty}\frac{d(\omega\ast \mu)(t)}{t^{\lambda+1-\alpha}}\\
 &=\Gamma(\lambda+1-\alpha)\int_0^{\infty}\int_0^{\infty}\frac{d\omega(u)}{(s+u)^{\lambda+1-\alpha}}\, d\mu(s),
 \end{align*}
 where $\omega$ is the measure in the Bernstein representation of the completely monotonic function $t^{-\lambda}f(t)$.
\end{cor}
{\it Proof.} We know from \cite[Corollary 2.1]{KP1} that the function $t^{-\lambda}f(t)$ is completely monotonic, and hence is of the form $t^{-\lambda}f(t)=\mathcal L(\omega)(t)$,
for some positive measure $\omega$ on $[0,\infty)$. This gives, using Lemma \ref{lemma:X},
\begin{align*}
\int_0^{\infty}f(t)t^{-\alpha}\mathcal L (\mu)(t)\, dt&=
\int_0^{\infty}t^{\lambda-\alpha}\mathcal L (\omega)(t)\mathcal L(\mu)(t)\, dt\\
&=
\Gamma(\lambda+1-\alpha)\int_0^{\infty}\frac{d (\omega\ast \mu)(s)}{s^{\lambda+1-\alpha}}.
\end{align*}
This proves the result.\hfill $\square$

Taking $\mu$ to be the point mass $\epsilon_x$ at $x$ the corollary yields the following.
\begin{cor}
 \label{lemma:F1} Let $f\in \bl$ and suppose that $\alpha<\lambda+1$. Then 
 $$
 \mathcal L(t^{-\alpha}f(t))(x)=\Gamma(\lambda+1-\alpha)\int_0^{\infty}\frac{d\omega(s)}{(x+s)^{\lambda+1-\alpha}},
 $$
 where  
 $\mathcal L(\omega)(t)=t^{-\lambda}f(t)$.
 In particular, the Laplace transform maps the class $t^{-\alpha}\bl$ into $\mathcal S_{\lambda+1-\alpha}$.

\end{cor}

Next, let us show some consequences of these results. The incomplete gamma function is, due to the relation \eqref{eq:int-rep-bl}, a fundamental building block in constructing functions in $\bl$. 
\begin{prop}
 \label{prop:incomplete-relation}
The following relation holds for any positive measure $\mu$, $\alpha<\lambda+1$ and $x>0$
\begin{align*}
\int_0^{\infty}\gamma(\lambda,xt)t^{-\alpha}\mathcal L (\mu)(t)\, dt&=\Gamma(\lambda+1-\alpha)\int_0^{\infty}\int_0^{x}\frac{u^{\lambda-1}du}{(s+u)^{\lambda+1-\alpha}}\, d\mu(s)\\
&=\Gamma(\lambda+1-\alpha)\int_0^{\infty}B\left(\lambda,1-\alpha; \frac{x}{x+s}\right)\, \frac{d\mu(s)}{s^{1-\alpha}}.
\end{align*}
\end{prop}
{\it Proof.}
Since the function $f(t)=\gamma(\lambda,xt)$ belongs to $\bl$ and 
$$
\gamma(\lambda,xt)t^{-\lambda}=\int_0^xe^{-tu}u^{\lambda-1}\, du,
$$
the corresponding measure in Corollary \ref{lemma:F1} is  $d\omega(u)=\mathbbm 1_{[0,x]}(u)u^{\lambda-1}du$. 
Applying now Corollary \ref{lemma:F3} 
we obtain 
\begin{equation*}
\int_0^{\infty}\gamma(\lambda,xt)t^{-\alpha}\mathcal L (\mu)(t)\, dt=\Gamma(\lambda+1-\alpha)\int_0^{\infty}\int_0^{x}\frac{u^{\lambda-1}du}{(s+u)^{\lambda+1-\alpha}}\, d\mu(s).
\end{equation*}
Since 
$$
\int_0^{x}\frac{u^{\lambda-1}du}{(s+u)^{\lambda+1-\alpha}}=\frac{1}{s^{1-\alpha}}B\left(\lambda,1-\alpha; \frac{x}{x+s}\right)
$$
(which follows by the change of variable $v=u/(u+s)$) for $\alpha<\lambda +1$, the proof is  complete.\hfill $\square$

Letting $\alpha=0$ and $\mu$ be the point mass at $s$ in Proposition \ref{prop:incomplete-relation} we record the following elementary relation (see also \cite[6.451.1]{GR}):
\begin{cor}
\label{cor:special-case}
$$
\int_0^{\infty}e^{-xt}\gamma(\lambda,ts)\, dt=\frac{\Gamma(\lambda)}{x}\frac{s^{\lambda}}{(x+s)^{\lambda}}.
$$
\end{cor}

Let us mention a couple of examples: Again letting $\mu$ be the point mass at $s$ in Proposition \ref{prop:incomplete-relation} we obtain the assertion in the motivating example mentioned in the introduction.
\begin{ex} 
\label{ex:fundamental}
For $s>0$ the function 
$$
x\mapsto B\left(\lambda,1-\alpha;\frac{x}{x+s}\right)
$$
belongs to $\mathcal T_{\lambda,\alpha}$.
\end{ex}
\begin{ex}
For a positive integer $n$ and $c>0$ the function 
$$
g_{c,n}(x)=(-1)^{n-1}\left( \log\left(\frac{x+c}{c}\right)+\sum_{k=1}^{n-1}\frac{(-1)^kx^k}{k}c^{-k}\right)
$$
belongs to $\mathcal T_{n,n}$. (See Proposition \ref{prop:Apage17} with $\mu=\epsilon_c$.)
Notice that 
$$
g_{c,n}(x)=x^n/c^n\, \int_0^{\infty}e^{-xt/c}E_n(t)\, dt,
$$
where $E_n$ is the generalized exponential integral (see \cite[8.19.24]{dlmf} and the end of Section \ref{sec:additional}).

Also the function
$$
x\mapsto \Gamma(\lambda)x^{\lambda}\int_0^{\infty}\int_{t_n}^{\infty}\cdots
\int_{t_2}^{\infty}\int_{t_1}^{\infty}\frac{ds}{(x+s)^{\lambda}s}dt_1\cdots dt_{n-1}d\mu(t_n)
$$
(where $\mu$ is any positive measure making this multiple integral converge) belongs to $\mathcal T_{\lambda,n}$. (See Proposition \ref{prop:handwritten} with $\beta=0$.)
\end{ex}

\begin{thm}
 \label{thm:A2} For a function $f:(0,\infty)\to [0,\infty)$ we have
 $$ f\in \tla \ \Leftrightarrow \ x^{1-\lambda}f'(x)\in \mathcal S_{\lambda+1-\alpha}.
 $$
\end{thm}
{\it Proof.} Assume that $f\in \tla$. In Definition \ref{defn:fundamental} the function  $t^{\alpha}\varphi(t)$ is the Laplace transform of a positive measure $\mu$, and using Proposition \ref{prop:incomplete-relation}, $f$ can be written in the form 
\begin{align*}
f(x)&=ax^{\lambda}+b+\int_0^{\infty}\gamma(\lambda,xt)t^{-\alpha}\mathcal L (\mu)(t)\, dt\\
&=ax^{\lambda}+b+\Gamma(\lambda+1-\alpha)\int_0^{\infty}\int_0^{x}\frac{u^{\lambda-1}du}{(s+u)^{\lambda+1-\alpha}}\, d\mu(s).
\end{align*}
This gives
$$
f'(x)=\lambda ax^{\lambda-1}+\Gamma(\lambda+1-\alpha)x^{\lambda-1} \int_0^{\infty}\frac{d\mu(s)}{(s+x)^{\lambda+1-\alpha}},
$$
showing that $x^{1-\lambda}f'(x)\in \mathcal S_{\lambda+1-\alpha}$.

Conversely, if $x^{1-\lambda}f'(x)\in \mathcal S_{\lambda+1-\alpha}$ then 
$$
f'(x)=x^{\lambda -1}\int_0^{\infty}e^{-xs}s^{\lambda-\alpha}\varphi(s)\, ds+cx^{\lambda-1},
$$
where $\varphi$ is completely monotonic and $c\geq 0$. Since $f$ is increasing, integration of this relation yields
\begin{align*}
f(x)-f(0+)&=\int_0^{\infty}\int_0^xt^{\lambda -1}e^{-ts}\, dts^{\lambda-\alpha}\varphi(s)ds+\frac{c}{\lambda}x^{\lambda}\\
&=\int_0^{\infty}\gamma(\lambda,sx)s^{-\alpha}\varphi(s)ds+\frac{c}{\lambda}x^{\lambda}.
\end{align*}
By definition, $s^{-\alpha}\varphi(s)$ is completely monotonic of order $\alpha$ and so $f\in \mathcal T_{\lambda,\alpha}$. This completes the proof.
\hfill $\square$
\begin{rem} Theorem \ref{thm:A2} yields
\begin{enumerate}[(a)]
 \item when $\alpha=1$: $f\in \mathcal T_{\lambda,1}\Leftrightarrow x^{1-\lambda}f'(x)\in \mathcal S_{\lambda}$ (see also \cite[Theorem 4.1]{KP1}); 
 \item when $\alpha=\lambda$: $f\in \mathcal T_{\lambda,\lambda} \Leftrightarrow x^{1-\lambda}f'(x)\in \mathcal S$;  
  \item when $\lambda=1$: $f\in \mathcal T_{1,\alpha} \Leftrightarrow f'\in \mathcal S_{2-\alpha}$. See also Remark \ref{rem:Apge21}.
\end{enumerate}
\end{rem}

\section{Representation via Hypergeometric functions}
In this section we prove that the $(\lambda,\alpha)$-Thorin-Bernstein functions admit an integral representation in terms of the hypergeometric function $\tFo$.
\begin{thm}
 \label{thm:A22}
 A function $f$ belongs to $\tla$ if and only if there are non-negative constants $a$ and $b$ and a positive measure $\mu$ such that
 $$
 f(x)=ax^{\lambda}+b+\frac{\Gamma(\lambda+1-\alpha)x^{\lambda}}{\lambda}\int_0^{\infty}\frac{\tFo(\alpha,1; \lambda+1; -x/s)}{(s+x)^{\lambda-\alpha}}\frac{d\mu(s)}{s}.
 $$
 In the affirmative case, $\mu$ is the measure such that 
 $t^{\alpha}\varphi(t)=\mathcal L(\mu)(t)$, $\varphi$ being the function in Definition \ref{defn:fundamental}.
\end{thm}
{\it Proof.} In the first relation in Proposition \ref{prop:incomplete-relation} we perform the change of variable $v=u/x$. This gives us
\begin{align*}
 f(x)-ax^{\lambda}-b&=\int_0^{\infty}\gamma(\lambda,xt)t^{-\alpha}\mathcal L (\mu)(t)\, dt\\&=\Gamma(\lambda+1-\alpha)x^{\lambda}\int_0^{\infty}\frac{1}{s^{\lambda+1-\alpha}}\int_0^{1}\frac{v^{\lambda-1}dv}{(1+vx/s)^{\lambda+1-\alpha}}\, d\mu(s).
 \end{align*}
Next, a combination of Euler's integral representation of the $\tFo$ and Euler's transformation (see \cite[Theorem 2.2.1 and Theorem 2.2.5]{AAR}) yields
\begin{align*}
 f(x)-ax^{\lambda}-b
 &=\frac{\Gamma(\lambda+1-\alpha)x^{\lambda}}{\lambda}\int_0^{\infty}\frac{\tFo(\lambda+1-\alpha,\lambda;\lambda+1;-x/s)}{s^{\lambda+1-\alpha}}\, d\mu(s)\\
 &=\frac{\Gamma(\lambda+1-\alpha)x^{\lambda}}{\lambda}\int_0^{\infty}\frac{\tFo(\alpha,1;\lambda+1;-x/s)}{s^{\lambda+1-\alpha}(1+x/s)^{\lambda-\alpha}}\, d\mu(s).
\end{align*}
This proves the asserted formula.\hfill $\square$

We record the following equivalent characterizations, obtained by using Pfaff's transformation (\cite[Theorem 2.2.5]{AAR}), Proposition \ref{prop:incomplete-relation} and the identity 
$$
B(c,d;z)=\frac{z^c}{c}\, \tFo\left(
c,1-d;c+1; z\right),
$$
see \cite[8.17.7]{dlmf}.
\begin{cor}
 \label{cor:A2}
 The following statements are equivalent for a function $f:(0,\infty)\to [0,\infty)$.
 \begin{enumerate}
  \item[(a)] $f\in \tla$,
  \item[(b)] $f$ can be represented as
  $$
 f(x)=ax^{\lambda}+b+\frac{\Gamma(\lambda+1-\alpha)x^{\lambda}}{\lambda}\int_0^{\infty}\frac{\tFo(\lambda, \alpha; \lambda+1; x/(x+s))}{(s+x)^{\lambda}}\frac{d\mu(s)}{s^{1-\alpha}},
 $$
 \item[(c)] $f$ can be represented as
  $$
 f(x)=ax^{\lambda}+b+\Gamma(\lambda+1-\alpha)\int_0^{\infty}B\left(\lambda,1-\alpha;\frac{x}{x+t}\right)\frac{d\mu(t)}{t^{1-\alpha}}.
 $$
 \end{enumerate} 
\end{cor}

\begin{rem} As was shown in Proposition \ref{prop:easy-pointwise}, if $f_n\in \tla$ converges pointwise to $f$ then $f$ also belongs to $\tla$. Letting $f_n$ correspond to the triple $(a_n,b_n,\mu_n)$ in (c) of Corollary \ref{cor:A2} and $f$ to $(a,b,\mu)$, then $\mu_n\to \mu$ vaguely on $(0,\infty)$.  
To see this we notice 
\begin{align*}
x^{1-\lambda}f_n'(x)&=\lambda a_n+\Gamma(\lambda+1-\alpha)\int_0^{\infty}\frac{d\mu_n(t)}{(x+t)^{\lambda+1-\alpha}}\quad \text{and}\\
x^{1-\lambda}f'(x)&=\lambda a+\Gamma(\lambda+1-\alpha)\int_0^{\infty}\frac{d\mu(t)}{(x+t)^{\lambda+1-\alpha}},
\end{align*}
from which it follows that $\mu_n\to \mu$ vaguely (see \cite[Theorem 2.2]{S}).\end{rem}
\begin{rem}
 \label{rem:Apage16} Suppose $\lambda=\alpha=1$. Since $\tFo(1,1;2;x)=-(1/x)\log (1-x)$ the function $f$ in Corollary \ref{cor:A2} takes the form 
 $$
 f(x)=ax+b+\int_0^{\infty}\log\left(\frac{x+t}{t}\right)\, d\mu(t),
 $$
 which is in accordance with the representation of ordinary Thorin-Bernstein functions. See \cite[Theorem 8.2]{S}. See also Proposition \ref{prop:Apage17}.
\end{rem}

\begin{rem}
 \label{rem:Apge21} 
 In the special case where $\lambda =1$ and $0<\alpha<1$ observe that for $-1<z<1$,
 $$
 z\, \tFo\left(1, \alpha; 2; z\right)=\int_0^z\frac{dt}{(1-t)^{\alpha}}.
 $$           
Substituting $z=x/(x+s)$ it follows that 
$$
\frac{x}{x+s}\, \tFo\left(1,\alpha; 2; \frac{x}{x+s}\right)=\frac{1}{1-\alpha}-\frac{1}{1-\alpha}\left(\frac{s}{x+s}\right)^{1-\alpha}.
$$
Therefore, Corollary \ref{cor:A2}(b) gives us that any $f\in \mathcal T_{1,\alpha}$ satisfies
$$
f'(x)=a+\Gamma(2-\alpha)\int_0^{\infty}\frac{d\mu(s)}{(x+s)^{2-\alpha}},
$$
which is \cite[Proposition 8.10]{S}. 
\end{rem}

\section{Integer values of the parameters}

We find in Proposition \ref{prop:Apage17} another form of the representation of functions from $\tla$ in the case where $\lambda=\alpha=n$, $n$ being a positive integer.  First, a lemma.
\begin{lemma}
 \label{lemma:Apage18}
 For $z\in (-1,1)\setminus \{0\}$ we have 
 $$
 \frac{\tFo(\lambda,\alpha;\lambda+1;z)}{\lambda}=\sum_{k=0}^{\infty}\frac{1}{\lambda+k}\frac{(\alpha)_k}{k!}z^k=
 \frac{1}{z^{\lambda}}\int_0^{z/(1-z)}u^{\lambda-1}(1+u)^{\alpha-\lambda-1}\, du.
 $$
 \end{lemma}
{\it Proof.} The first equality holds by definition. Furthermore, 
 \begin{align*}
 \sum_{k=0}^{\infty}\frac{1}{\lambda+k}\frac{(\alpha)_k}{k!}z^k
 &=
 \frac{1}{z^{\lambda}}\sum_{k=0}^{\infty}\frac{1}{\lambda+k}\frac{(\alpha)_k}{k!}z^{k+\lambda}\\
 &=
 \frac{1}{z^{\lambda}}\int_0^zt^{\lambda-1}\sum_{k=0}^{\infty}\frac{(\alpha)_k}{k!}t^{k}\, dt\\
 &=
 \frac{1}{z^{\lambda}}\int_0^zt^{\lambda-1}\frac{1}{(1-t)^{\alpha}}\, dt.
 \end{align*}
 The change of variable $u=t/(1-t)$ transforms this expression into
 $$
 \frac{1}{z^{\lambda}}\int_0^{z/(1-z)}u^{\lambda-1}(1+u)^{\alpha-\lambda-1}\, du,
 $$
 and this proves the lemma.\eop

 \begin{prop}
 \label{prop:Apage17}
 Suppose that $n$ is a positive integer. A function $f$ belongs to $\mathcal T_{n,n}$ if and only if 
 $$
 f(x)=ax^{n}+b+(-1)^{n-1}\int_0^{\infty}\left( \log\left(\frac{x+s}{s}\right)+\sum_{k=1}^{n-1}\frac{(-1)^kx^k}{ks^k}\right)\, s^{n-1}\, d\mu(s).
 $$
\end{prop}
{\it Proof.} Corollary \ref{cor:A2} gives us 
$$
f(x)=ax^{n}+b+\frac{x^{n}}{n}\int_0^{\infty}\frac{\tFo(n, n; n+1; x/(x+s))}{(s+x)^{n}}\frac{d\mu(s)}{s^{1-n}},
$$
and Lemma \ref{lemma:Apage18} yields (for $x\neq 0$)
$$
\frac{\tFo(n, n; n+1; x/(x+s))}{n}=\frac{(s+x)^{n}}{x^n}\int_0^{x/s}\frac{u^{n-1}}{1+u}\, du.
$$
The integral in the right-hand side of this relation can be rewritten using Taylor's formula with integral remainder:
$$
\int_0^{x/s}\frac{u^{n-1}}{1+u}\, du=(-1)^{n-1}\left(\log (1+x/s)+\sum_{k=1}^{n-1}\frac{(-x/s)^k}{k}\right).
$$
In this way the proposition is proved.\eop

\begin{rem}
In the case where $\alpha=1$ and $\lambda$ is a positive rational number, the integrand in the representation in Corollary \ref{cor:A2}(b) can be summed. 
Indeed, writing $\lambda=n+p/q$, where $n, p$ and $q$ are non-negative integers and $p<q$ we have
\begin{align*}
   \sum_{k=0}^{\infty}\frac{1}{k+\lambda}\left(\frac{x}{x+t}\right)^{k+\lambda}&=
   \sum_{k=0}^{\infty}\frac{1}{k+p/q}\left(\frac{x}{x+t}\right)^{k+p/q}\\
   &\phantom{=}-
   \sum_{k=0}^{n-1}\frac{1}{k+p/q}\left(\frac{x}{x+t}\right)^{k+p/q}.
\end{align*}
The infinite sum on the right hand side can be rewritten as
$$
\sum_{k=0}^{\infty}\frac{1}{k+p/q}\left(\frac{x}{x+t}\right)^{k+p/q}=-\sum_{k=0}^{q-1}e^{-2\pi i kp/q}\log \left(1-\left(\frac{x}{x+t}\right)^{1/q}e^{2\pi i k/q}\right).
$$
This formula can be found in \cite{Prudnikov}. We indicate how to prove a more general formula below. The representation thus takes the form 
\begin{align*}
f(x)=ax^{n+p/q}+b&+\\
\Gamma(n+p/q)\int_0^{\infty}\Bigg\{ &-\sum_{k=0}^{q-1}e^{-2k\pi ip/q}\log\left(1-\left(\frac{x}{x+t}\right)^{1/q}e^{2k\pi i/q}\right)\\
&-\left(\frac{x}{x+t}\right)^{p/q}\sum_{k=0}^{n-1}\frac{1}{k+p/q}\left(\frac{x}{x+t}\right)^k\Bigg\} \, d\mu(t).
\end{align*}
In particular, for $\lambda=n+1/2$, using  $2\tanh ^{-1}(s)=\log ((1+s)/(1-s))$ for $|s|<1$,
\begin{align*}
f(x)=ax^{n+1/2}+b&+
2\Gamma(n+1/2)\int_0^{\infty}\left\{ \tanh ^{-1}\left(\left(\frac{x}{x+t}\right)^{1/2}\right)\right.\\
&-\left. \sum_{k=0}^{n-1}\frac{1}{2k+1}\left(\frac{x}{x+t}\right)^{k+1/2}\right\} \, d\mu(t).
\end{align*}

The summation formula alluded to: 
Suppose that $g(z)=\sum_{n=1}^{\infty}a_nz^n$ converges for $|z|<1$. Then, interchanging the order of summation and using elementary properties of the roots of unity it follows that
$$
\sum_{k=0}^{q-1}e^{-2\pi i kp/q}g\left(z^{1/q}e^{2\pi i k /q}\right)=
qz^{p/q}\sum_{k=0}^{\infty}a_{p+kq}z^k.
$$
The formula needed above corresponds to $a_n=1/n$.
\end{rem}
In the next proposition we examine in more detail the representation of functions in $\tla$ taking into account both the integer part and the fractional part of the order $\alpha$ of complete monotonicity of $\varphi$.

\begin{prop}
\label{prop:handwritten}
 Suppose that $\varphi$ is a completely monotonic function of order $\alpha=n+\beta$, $n\in \mathbb N$, $\beta\in [0,1)$. There exists a positive measure $\mu$ such that, with
$$
\xi_n(t)=\frac{1}{(n-1)!}\int_0^t(t-u)^{n-1}\,d\mu(u),
$$
we have
\begin{align*}
\lefteqn{\int_0^{\infty}\gamma(\lambda,xt)\varphi(t)\, dt}\\
&=\Gamma(\lambda-\beta+1)\int_0^{\infty}B\left(\lambda,1-\beta; \frac{x}{x+t}\right)\frac{\xi_n(t)}{t^{1-\beta}}\,dt\\
&=\Gamma(\lambda-\beta+1)\int_0^{\infty}\int_{t_n}^{\infty}\cdots
\int_{t_1}^{\infty}B\left(\lambda,1-\beta; \frac{x}{x+s}\right)\frac{ds}{s^{1-\beta}}dt_1\cdots d\mu(t_n),
\end{align*}
\end{prop}
{\it Proof.} According to \cite{KP0} we may write $\varphi$ as $\varphi(t)=t^{-\beta}\mathcal L(\xi_n)(t)$ where $\xi_n$ is the fractional integral of positive integer order 
$$
\xi_{n}(s)=\frac{1}{(n-1)!}\int_0^s(s-u)^{n-1}\, d\mu(u).
$$
The first equality now follows from Proposition \ref{prop:incomplete-relation}.

As it is well known, we have 
$$
\frac{1}{(n-1)!}\int_0^t(t-u)^{n-1}\, d\mu(u)=\int_0^t\int_0^{t_1}\cdots\int_0^{t_{n-1}}d\mu(t_n)dt_{n-1}\cdots dt_1
$$
and this gives the second equality.  \eop

\begin{rem}
When $n=0$ the contents of Proposition \ref{prop:handwritten} are described in Proposition \ref{prop:incomplete-relation}.
\end{rem}

Let us complete this section by relating our results with earlier results for generalized complete Bernstein functions and generalized Thorin-Bernstein functions presented in \cite{KP1}.

\begin{rem}
Proposition \ref{prop:handwritten} extends two results from \cite{KP1}:
\begin{enumerate}
\item[(a)] A function $f$ belongs to $\mathcal T_{\lambda,0}$ if and only if
$$
f(x)=ax^\lambda + b+ 
\Gamma(\lambda)\,\int_0^{\infty}\left(\frac{x}{x+s}\right)^{\lambda}\,\frac{d\mu(s)}{s}.
$$
(See also \cite[Proposition 3.12]{KP1}.)
\item[(b)] A function $f$ belongs to $\mathcal T_{\lambda,1}$ if and only if
$$
f(x)=ax^\lambda + b+ 
\Gamma(\lambda)\,\int_0^{\infty}\left(\frac{x}{x+s}\right)^{\lambda}\,\frac{h(s)}{s}\, ds,
$$
where $h=\xi_1$ is non-negative and increasing. (See also \cite[Proposition 4.1]{KP1}.)
\end{enumerate}
\end{rem}

\section{Additional results and comments}
\label{sec:additional}
For non-negative $\alpha$, it is easily seen that $(-1)^n\varphi^{(n)}$ is completely monotonic of order $\alpha$ if $\varphi$ is completely monotonic of order $\alpha$.  
Indeed, this follows by writing $\varphi(t)=t^{-\alpha}\mathcal L(\sigma)(t)$ and using Leibniz' formula:
$$
t^{\alpha}(-1)^n\varphi^{(n)}(t)
=\sum_{k=0}^n\binom{n}{k}(\alpha)_k\frac{1}{t^k}\mathcal L(s^{n-k}d\sigma(s))(t).
$$
Denoting by $m_r$  the measure on $(0,\infty)$ having density $s^{r-1}/\Gamma(r)$ w.r.t.\ Lebesgue measure (for $r>0$) we notice that $\mathcal L(m_r)(t)=t^{-r}$.
Thus  the relation above can be written as
$$
t^{\alpha}(-1)^n\varphi^{(n)}(t)=\mathcal L \left(s^nd\sigma(s)+\sum_{k=1}^n\binom{n}{k}(\alpha)_k (m_{k}\ast s^{n-k}d\sigma(s))\right)(t).
$$

\begin{defn}
For $\alpha\in [0,\lambda)$ and a non-negative integer $n$, $\tla^{(n)}$ is the subclass of $\tla$ consisting of the functions of the form
\begin{equation}
\label{eq:Apage1}
f(x)=a x^{\lambda}+b+
\int_0^{\infty}\gamma(\lambda,xt)(-1)^n\varphi^{(n)}(t)\, dt,
\end{equation}
where $a$ and $b$ are non-negative numbers, and 
$\varphi$ is a completely monotonic function of order $\alpha$.
\end{defn}
(It should be noted that the condition $\alpha<\lambda$ appears because that e.g.\  $-\varphi'(t)$ dominates $t^{-\alpha -1}$ for $t$ near $0$.)

\begin{prop} 
Suppose that $n\geq 1$. If $f$ belongs to $\tla^{(n)}$ with the representation \eqref{eq:Apage1} then $f(x)/x^{\lambda}$ is in $\mathcal S_{\lambda-\alpha}$. It has the representation
\begin{equation}
\label{eq:Apage1+1}
 f(x)/x^\lambda=a +bx^{-\lambda}+\Gamma(\lambda-\alpha)\int_0^{\infty}\frac{d\sigma_n(t)}{(x+t)^{\lambda-\alpha}},
\end{equation}
where 
\begin{equation}
    \label{eq:sigma_n}
\sigma_n=s^{n-1}d\sigma(s)+\sum_{k=1}^{n-1}\binom{n-1}{k}(\alpha)_k\, m_k\ast (s^{n-1-k}d\sigma(s)),
\end{equation}
$\sigma$ being the positive measure such that $t^{\alpha}\varphi(t)=\mathcal L(\sigma)(t)$.

Conversely, if $a,b\geq 0$ and $\sigma$ is any positive measure for which the Laplace transform converges then defining the measure $\sigma_n$ by \eqref{eq:sigma_n}, $f$ given by \eqref{eq:Apage1+1}  belongs to $\tla^{(n)}$.
The corresponding function $\varphi$ in \eqref{eq:Apage1} is given by $\varphi(t)=t^{-\alpha}\mathcal L(\sigma)(t)$ and the measure $\sigma_n$ is related to $\varphi$ in the following way:
$(-1)^{n-1}\varphi^{(n-1)}(t)=t^{-\alpha}\mathcal L(\sigma_{n})(t)$.

\end{prop}
{\it Proof.} When $\varphi$ is a completely monotonic function of order $\alpha\geq 0$ we write $t^{\alpha}\varphi(t)=\mathcal L(\sigma)(t)$, $\varphi(t)=\mathcal L(\mu)(t)$,  and notice that $(-1)^n\varphi^{(n)}(t)=\mathcal L(s^nd\mu(s))(t)$. 
This gives, using Lemma \ref{lemma:F2} and Corollary \ref{cor:special-case},
\begin{align*}
\int_0^{\infty}\gamma(\lambda,xt)(-1)^n\varphi^{(n)}(t)\, dt&=
\int_0^{\infty}\gamma(\lambda,xt)\mathcal L(s^nd\mu(s))(t)\, dt\\
&=\Gamma(\lambda)x^{\lambda}\int_0^{\infty}\frac{s^{n-1}}{(s+x)^{\lambda}}\, d\mu(s)\\
&=x^{\lambda}\int_0^{\infty}e^{-xt}t^{\lambda-1}\mathcal L(s^{n-1}d\mu(s))(t)\, dt\\
&=x^{\lambda}\int_0^{\infty}e^{-xt}t^{\lambda-1}(-1)^{n-1}\varphi^{(n-1)}(t)\, dt\\
&=x^{\lambda}\int_0^{\infty}e^{-xt}t^{\lambda-1-\alpha}\mathcal L(\sigma_{n})(t)\, dt\\
&=\Gamma(\lambda-\alpha)x^{\lambda}\int_0^{\infty}\frac{d\sigma_{n}(t)}{(x+t)^{\lambda-\alpha}},
\end{align*}
where $\sigma_{n}$ is the positive measure given in \eqref{eq:sigma_n} and thus  satisfies 
$$(-1)^{n-1}\varphi^{(n-1)}(t)=t^{-\alpha}\mathcal L(\sigma_{n})(t).
$$
The proof of the other direction follows by retracing these steps. This concludes the proof of the proposition.\eop

A $C^{\infty}$-function $f:(0,\infty)\to (0,\infty)$ is said to be logarithmically completely monotonic if $-(\log f(x))'=-f'(x)/f(x)$ is completely monotonic. A measure or function on the positive half line is called infinitely divisible if its Laplace transform is a logarithmically completely monotonic function. For more details see \cite{BKP}.

\begin{cor}
Let $\sigma$ and $\sigma_n$ be as in the proposition above. Then the function $t^{\lambda-\alpha-1}\mathcal L(\sigma_n)(t)$ is infinitely divisible if $1<\lambda-\alpha\leq 2$.
\end{cor}
{\it Proof.} It follows from the relation 
$$
\int_0^{\infty}e^{-xt}t^{\lambda-1-\alpha}\mathcal L(\sigma_{n})(t)\, dt=\Gamma(\lambda-\alpha)\int_0^{\infty}\frac{d\sigma_{n}(t)}{(x+t)^{\lambda-\alpha}}
$$ 
in the proof above that the Laplace transform of $t^{\lambda-\alpha-1}\mathcal L(\sigma_n)(t)$ belongs to $\mathcal S_{\lambda-\alpha}$ and thus to $\mathcal S_2$. Therefore it is logarithmically completely monotonic by a result of Kristiansen. See \cite{K} and \cite{BKP}. \eop

Let us end this section with a few additional observations on infinite divisibility.

 First of all, if $f\in C_{\alpha}$ for some $\alpha\geq -1$ then $f$ is infinitely divisible. This is well-known for $\alpha\geq 0$ and for  $\alpha\in [-1,0)$ it follows from the relation
 $$
 \mathcal L(f)(x)=\int_0^{\infty}e^{-xt}t^{1-\alpha-1}f(t)t^{\alpha}\, dt
 $$
 and Kristiansen's theorem.

The next results deal with products of completely monotonic functions of given real orders. 
The following is a consequence of Lemma \ref{lemma:X}.
\begin{prop}
 \label{prop:Gpage1} If $f\in \mathcal C_{\alpha}$, $g\in \mathcal C_{\beta}$ with $\alpha+\beta<1$ and 
 $$
 f(t)=t^{-\alpha}\mathcal L(\mu)(t),\quad g(t)=t^{-\beta}\mathcal L(\nu)(t)
 $$
 then $\mathcal L(fg)\in \mathcal S_{1-\alpha-\beta}$ and 
 $$
 \int_0^{\infty}e^{-xt}f(t)g(t)\, dt=\Gamma(1-\alpha-\beta)\int_0^{\infty}\frac{d(\mu\ast\nu)(t)}{(x+t)^{1-\alpha-\beta}}.
 $$
\end{prop}

\begin{cor}
 \label{cor:Gpage4}
 If $f\in \mathcal C_{\alpha}$, $g\in \mathcal C_{\beta}$ and $-1\leq \alpha+\beta<1$, then the function $fg$ is infinitely divisible. 
\end{cor}
{\it Proof.} From Proposition \ref{prop:Gpage1} we see that $\mathcal L(fg)$ belongs to $\mathcal S_{1-\alpha-\beta}$, which is contained in $\mathcal S_2$, and Kristiansen's result can thus be used. \eop

The next corollary deals also with the distribution function 
$$
F_{\lambda}(t)=1-\lambda t^{\lambda}e^t\Gamma(-\lambda,t)$$
of a so-called randomized Lomax distribution, defined in terms of the complementary incomplete Gamma function,
$$
\Gamma(\lambda,x)=\int_x^{\infty}e^{-u}u^{\lambda -1}\, du.
$$ See \cite[Example 3.11]{KP1}, where it was shown that $F_{\lambda}\in \bl$. 
A standard computation shows that
$$
F_{\lambda}(x)=\frac{1}{\Gamma(\lambda)}\int_0^{\infty}\gamma(\lambda,xt)\frac{dt}{(1+t)^2},
$$
from which it is immediate that $F_{\lambda}$ even belongs to $\mathcal T_{\lambda,0}$.
\begin{cor}
 \label{cor:Fpage7}
 Suppose that $\lambda>0$ and $\lambda-1\leq\alpha< \lambda+1$. If $f\in\bl$ then $t^{-\alpha}f(t)$ is infinitely divisible. In particular, $t^{-\alpha}\gamma(\lambda,xt)$ (for fixed $x$) and $t^{-\alpha}F_{\lambda}(t)$ are infinitely divisible.
\end{cor}
{\it Proof.} Of course $t^{-\alpha}$ belongs to $\mathcal C_{\alpha}$ and $f$ belongs to $\mathcal C_{-\lambda}$. The result now follows from Corollary \ref{cor:Gpage4}.\eop

Next we investigate the image of the subclasses $\tla$  under the Laplace transform. 
\begin{prop} 
\label{prop:dd}
    If $f\in\tla$ then $x^{\alpha}\mathcal L(f)(x)-f(0+)x^{\alpha-1}$ belongs to the class $\mathcal S_{\lambda+1-\alpha}$. 
\end{prop}
\begin{rem}
If $\alpha\leq 1$ then $x^{\alpha}\mathcal L(f)(x)$ in the proposition above belongs to $\mathcal S_{\lambda+1-\alpha}$. In particular, $\mathcal L(f)$ belongs to $\mathcal S_{\lambda+1}$ for $f\in \mathcal T_{\lambda,0}$. Hence, if $\lambda\leq 1$ then any function  $f\in \mathcal T_{\lambda,0}$ is infinitely divisible.
\end{rem}
{\it Proof of Proposition \ref{prop:dd}.} 
Let $f\in \tla$ and write its representation as
$$
f(t)=at^{\lambda}+b+\int_0^{\infty}\gamma(\lambda,ts)s^{-\alpha}\mathcal L(\mu)(s)\, ds
$$
with $b=f(0+)$ (see \cite[Proposition 2.2]{KP1}) and some positive measure $\mu$. This gives, using Corollary \ref{cor:special-case} and Fubini's theorem,
$$
\mathcal L(f)(x)=\frac{a\Gamma(\lambda+1)}{x^{\lambda+1}}+\frac{b}{x}+\frac{\Gamma(\lambda)}{x}\int_0^{\infty}s^{-\alpha}\mathcal L(\mu)(s)\frac{s^{\lambda}}{(x+s)^{\lambda}}\, ds.
$$
Notice that $\varphi(s)=s^{\lambda}/(x+s)^{\lambda}$
is in $\bl$ and that
$$
s^{-\lambda}\varphi(s)=\frac{1}{\Gamma(\lambda)}\int_0^{\infty}e^{-su}e^{-xu}u^{\lambda-1}\, du.
$$
According to Corollary \ref{lemma:F3} we obtain
$$
\int_0^{\infty}s^{-\alpha}\mathcal L(\mu)(s)\varphi(s)\, ds
=
\frac{\Gamma(\lambda+1-\alpha)}{\Gamma(\lambda)}\int_0^{\infty}\int_0^{\infty}\frac{u^{\lambda-1}e^{-xu}}{(u+s)^{\lambda+1-\alpha}}\, du \, d\mu(s),
$$
and thus
\begin{align*}
    x^{\alpha}\mathcal L(f)(x)&=\frac{a\Gamma(\lambda+1)}{x^{\lambda+1-\alpha}}+bx^{\alpha-1}\\
    &\phantom{=}+\Gamma(\lambda+1-\alpha)x^{\alpha-1}\int_0^{\infty}\int_0^{\infty}\frac{u^{\lambda-1}e^{-xu}}{(u+s)^{\lambda+1-\alpha}}\, du \, d\mu(s).
\end{align*}
In this last integral we make a change of variable $v=(x/s)u$ and get by Fubini's theorem 
\begin{align*}
    x^{\alpha}\mathcal L(f)(x)&=\frac{a\Gamma(\lambda+1)}{x^{\lambda+1-\alpha}}+bx^{\alpha-1}\\
    &\phantom{=}+\Gamma(\lambda+1-\alpha)\int_0^{\infty}\int_0^{\infty}s^{\alpha-1}e^{-vs}\, d\mu(s)\frac{v^{\lambda-1}dv}{(x+v)^{\lambda+1-\alpha}}.
\end{align*}
It is now clear that $x^{\alpha}\mathcal L(f)(x)-bx^{\alpha-1}$ belongs to $\mathcal S_{\lambda+1-\alpha}$.
\eop

We remark that in general not all functions from $\mathcal S_{\lambda+1-\alpha}$ are of the form $x^{\alpha}\mathcal L(f)(x)$, for some $f\in \tla$. For example, since $\mathcal T_{\lambda,1}\subsetneq \bl$ it follows that 
$\Phi(\mathcal T_{\lambda,1})\subsetneq \Phi(\bl)=\mathcal S_{\lambda}$. Below an example is given when $\lambda=\alpha\in \{1,2,3,\ldots\}$.
\begin{ex}
Let $n\geq 1$. To see that $\{x^n\mathcal L(g)(x)\, |\, g\in \mathcal T_{n,n}\}$ is not all of $\mathcal S$ let $\varphi$ be a completely monotonic function and suppose that there is $g\in \mathcal T_{n,n}$ such that $x^n\mathcal L(g)(x)=\mathcal L(\varphi)(x)$. This entails 
$$
g(x)=\int_0^x(x-s)^{n-1}\varphi(s)\, ds.
$$
According to Theorem \ref{thm:A2}, $g$ belongs to $\mathcal T_{n,n}$ if and only if $x^{1-n}g'(x)$ belongs to $\mathcal S$. We obtain
$$
x^{1-n}g'(x)=(n-1)\int_0^1(1-t)^{n-2}\varphi(xt)\, dt,
$$
which is not in general a Stieltjes function. As an example, take $\varphi(s)=e^{-s}$. In this case, $x^{1-n}g'(x)$ is the Laplace transform of the non-completely monotonic function $(n-1)(1-t)^{n-2} \mathbbm 1_{(0,1)}(t)$.
\end{ex}

Let us conclude the paper by giving some examples related to the generalized exponential integral, $E_p$, defined for $p>0$, as 
$$
E_p(x)=x^{p-1}\int_x^{\infty}\frac{e^{-t}}{t^p}\, dt
$$
We write
 $$
 E_p(x)=\int_1^{\infty}e^{-xt}\frac{dt}{t^{p}}
 $$
 and this yields that $E_p$ is a completely monotonic function for any $p>0$.
 We observe that 
$$
0\leq tE_p(t)\leq \int_t^{\infty}e^{-s}\, ds=e^{-t},
$$
so $E_p(t)$ decays exponentially for $t$ tending to $\infty$ and is bounded by $1/t$ as $t$ tends to $0$.

 \begin{ex}
     We have, using Proposition \ref{prop:incomplete-relation},
     $$
     \int_0^{\infty}\gamma(\lambda,xt)E_p(t)\, dt=\Gamma(\lambda)x^{\lambda}\int_1^{\infty}\frac{1}{(x+s)^{\lambda}}\frac{ds}{s^{p+1}}.$$
     The change of variable $t=x/(x+s)$ transforms the right-hand side into
     \begin{equation}
     \label{eq:5.10}
        \frac{\Gamma(\lambda)}{x^p}B(\lambda+p,-p; x/(x+1)),
     \end{equation}
     which is thus an example of a function in  $\mathcal T_{\lambda,0}$, for any $p>0$. Notice that it is infinitely divisible for $\lambda\leq 1$.
 \end{ex}

 \begin{ex}
  According to \cite[8.19.4]{dlmf} we have
$$
x^{1-p}E_p(x)=\frac{1}{\Gamma(p)}\int_0^{\infty}\frac{e^{-x(t+1)}t^{p-1}}{1+t}\, dt,
$$
 and it follows that $E_p\in \mathcal C_{1-p}$ with $E_p(x)=x^{p-1}\mathcal L(h_p)(x)$, where
$$
h_p(t)=\mathbbm 1_{(1,\infty)}(t)\frac{(t-1)^{p-1}}{\Gamma(p)t}.
$$
This gives, applying Proposition \ref{prop:incomplete-relation} followed by the change of variable $u=1-1/t$
\begin{align*}
\int_0^{\infty}\gamma(\lambda,xt)E_p(t)\, dt
&=
\Gamma(\lambda+p)\int_0^{\infty}B\left(\lambda,p; \frac{x}{x+t}\right)\frac{h_p(t)}{t^{p}}\, dt\\
&=\frac{\Gamma(\lambda+p)}{\Gamma(p)}\int_0^{1}B\left(\lambda,p; \frac{x(1-u)}{1+x(1-u)}\right)u^{p-1}\, du.
\end{align*}
Moreover, using the relation \eqref{eq:5.10} 
we obtain the identity 
$$
\frac{\Gamma(\lambda+p)}{\Gamma(p)}\int_0^{1}B\left(\lambda,p; \frac{x(1-u)}{1+x(1-u)}\right)u^{p-1}\, du=
  \frac{\Gamma(\lambda)}{x^p}B(\lambda+p,-p; x/(x+1))
$$
and this function belongs to $\mathcal T_{\lambda,1-p}$.
 \end{ex}
 \begin{ex}
     The function
$$
C_p(x)=\int_0^{\infty}\gamma(\lambda,xt)t^{1-p}E_p(t)\, dt
$$
is defined when $p<\lambda +1$. 
Furthermore, $C_p(x)$ belongs to $\mathcal T_{\lambda,0}$
and it has the representation 
$$
C_p(x)=\Gamma(\lambda)x^{\lambda}\int_0^{\infty}\frac{1}{(x+s)^{\lambda}}\frac{h_p(s)}{s}\, ds=\frac{\Gamma(\lambda)x^{\lambda}}{\Gamma(p)}\int_1^{\infty}\frac{(s-1)^{p-1}}{(x+s)^{\lambda}s^2}\, ds.
$$
The change of variable $t=(s-1)/s$, Euler's  integral representation of $\tFo$ and Pfaff's transformation yield
\begin{align*}
C_p(x)&=
\frac{\Gamma(\lambda)}{\Gamma(p)}x^{\lambda}\int_0^1t^{\lambda-p+1}(1-t)^{p-1}(1+xt)^{-\lambda}\, dt\\
 &=\frac{\Gamma(\lambda-p+2)}{\lambda(\lambda+1)}x^{\lambda}
  \tFo\left(
  \lambda,\lambda-p+2;
  \lambda+2
;-x\right)\\
 &=\frac{\Gamma(\lambda-p+2)}{\lambda(\lambda+1)}\left(\frac{x}{x+1}\right)^{\lambda}\tFo\left(
  \lambda, p;
  \lambda+2
;\frac{x}{x+1}\right).
\end{align*}
 \end{ex}

\noindent
Stamatis Koumandos\\
Department of Mathematics and Statistics\\
The University of Cyprus\\
P. O. Box 20537\\
1678 Nicosia, Cyprus\\
email: skoumand@ucy.ac.cy
\medskip

\noindent
Henrik Laurberg Pedersen (corresponding author)\\
Department of Mathematical Sciences\\
University of Copenhagen\\
Universitetsparken 5\\
DK-2100, Denmark\\
email: henrikp@math.ku.dk

\end{document}